\documentclass[a4paper,onecolumn]{article}

\usepackage{graphicx}
\usepackage{amssymb}
\usepackage{hyperref}
\begin{document}

{\setlength\arraycolsep{1pt}

\title{\textbf{Measurement of areas on a sphere using\\
Fibonacci and latitude--longitude lattices}}

\author{\textbf{\'Alvaro Gonz\'alez}\\
\\
\small{Departamento de Ciencias de la Tierra}\\
\small{Universidad de Zaragoza}\\
\small{C. Pedro Cerbuna, 12}\\
\small{50009 Zaragoza(Spain)}\\
\small{\href{mailto:Alvaro.Gonzalez@unizar.es}
{Alvaro.Gonzalez@unizar.es}}}

\date{\textbf{\textit{Mathematical Geosciences}, in press}\\
\small{\href{http://dx.doi.org/10.1007/s11004-009-9257-x}{http://dx.doi.org/10.1007/s11004-009-9257-x}}}

\maketitle

\begin{abstract}

The area of a spherical region can be easily measured by
considering which sampling points of a lattice are located inside
or outside the region. This point-counting technique is frequently
used for measuring the Earth coverage of satellite constellations,
employing a latitude--longitude lattice. This paper analyzes the
numerical errors of such measurements, and shows that they could
be greatly reduced if the Fibonacci lattice were used instead. The
latter is a mathematical idealization of natural patterns with
optimal packing, where the area represented by each point is
almost identical. Using the Fibonacci lattice would reduce the
root mean squared error by at least 40\%. If, as is commonly the
case, around a million lattice points are used, the maximum error
would be an order of magnitude smaller.

\end{abstract}

\section{Introduction}

The area of a region is easy to measure, without explicitly
considering its boundaries, by determining which points of a
lattice are inside or outside the region. This point-counting
method is commonly applied to estimating areas on a plane
\cite{995-Bardsley,996-Jarai,998-Howarth,997-Gundersen,Baddeley2004}.
A related issue is how to approximate the region boundaries from
this sampling \cite{972-Barclay}.

Point-counting on the sphere is commonly used for estimating the
Earth coverage of satellite constellations
\cite{939-Kantsiper,938-Feng}, the fraction of the Earth's surface
efficiently seen by one or more satellites. In the simplest case,
each satellite covers a circular region (cap), so the
constellation covers a complex set of isolated and/or overlapping
caps \cite{939-Kantsiper}. Similarly, some maps designed for
global earthquake forecasting depict earthquake-prone regions as
complex sets of up to tens of thousands of caps
\cite{Kossobokov2003,943-Kafka,GonzalezStatsei6}. To assess these
forecasts it is necessary to measure the fraction of the Earth's
area covered by these regions \cite{943-Kafka}.

Analytical solutions exist if the spherical region has a known,
regular or polygonal shape
\cite{971-Kimerling,940-Bevis,917-Sjoberg}. The area of a set of
caps on the sphere has a complex analytical solution
\cite{939-Kantsiper}, which unfortunately does not indicate
whether any particular location on the surface of the sphere is
covered.

The numerical error of point counting should ideally decrease
rapidly as the lattice density increases. Numerous works deal with
errors on the plane
\cite{995-Bardsley,996-Jarai,998-Howarth,997-Gundersen,Baddeley2004}.
Some particular cases using latitude--longitude lattices were
analyzed elsewhere \cite{939-Kantsiper}. In automated counting,
the computation time is directly proportional to the number of
lattice points. Satellite constellation coverage
\cite{937-Ochieng} or the areas marked on rapidly-updated
earthquake forecasting maps \cite{GonzalezStatsei6} need
continuous monitoring, implying a trade-off between accuracy and
computational load. Thus it is important to find lattices able to
measure areas as efficiently as possible.

On the plane, the regular hexagonal lattice provides optimal
sampling \cite{Conway1998}. On the sphere, it is impossible to
arrange regularly more than 20 points (the vertices of a
dodecahedron), and the optimum configuration of a large number of
points is problem-specific
\cite{916-Saff,Conway1998,948-Williamson,Gregory2008}. For optimal
point-counting, the area represented by every point should be
almost the same.

\begin{figure}[t]
\begin{center}
\includegraphics[width=\textwidth]{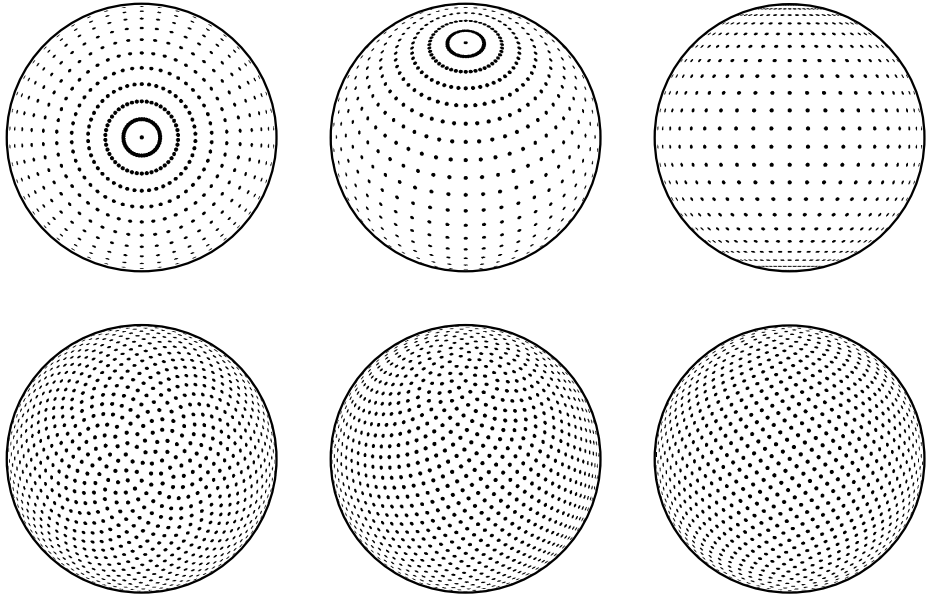}
\end{center}
\textbf{\caption{Latitude--longitude lattice (\textit{top}) and
Fibonacci lattice (\textit{bottom})\textnormal{, with 1014 and
1001 points, respectively. Orthographic projections, centred at
the pole (\textit{left}), latitude $45^\circ$ (\textit{middle})
and equator (\textit{right}). In the Fibonacci lattice, the points
are much more evenly spaced, and the axial anisotropy is much
smaller.}}\label{fig:SphereLattices}}
\end{figure}

Traditionally, the latitude--longitude lattice is used for
measuring Earth coverage
\cite{939-Kantsiper,937-Ochieng,938-Feng}. However, it is very
inhomogeneous (Fig.~\ref{fig:SphereLattices}), requiring
non-uniform weighting of the point contributions. Also, its number
of points is restricted by geometrical constraints.

The Fibonacci lattice is a particularly appealing alternative
\cite{Dixon1987,Dixon1989,Dixon1992,988-Fowler,974-Svergun,962-Kozin,
Swinbank1999,Swinbank2006a,902-Swinbank,978-Winfield,977-Nye,901-Hannay,970-Purser,
967-Purser}. Being easy to construct, it can have any odd number
of points \cite{902-Swinbank}, and these are evenly distributed
(Fig.~\ref{fig:SphereLattices}) with each point representing
almost the same area. For the numerical integration of continuous
functions on a sphere, it has distinct advantages over other
lattices \cite{901-Hannay,970-Purser}.

This paper demonstrates that for measuring the areas of spherical
caps, the Fibonacci lattice is much more efficient than its
latitude--longitude counterpart. After describing both lattices
(Sects.~\ref{sec:LatLonLattice} and \ref{sec:FiboLattice}), it is
explained how to use them for measuring cap areas
(Sect.~\ref{sec:AreaMeasure}) and how the error of this
measurement is assessed (Sect.~\ref{sec:ErrorAssess}). The error
results obtained from an extensive Monte Carlo simulation are
described, and to some extent explained analytically
(Sect.~\ref{sec:ErrResults}). The conclusions are set out in the
final section.

\section{Latitude--Longitude Lattice}\label{sec:LatLonLattice}

The latitude--longitude lattice is the set of points located at
the intersections of a grid of meridians and parallels, separated
by equal angles of latitude and longitude
(Fig.~\ref{fig:SphereLattices}). This is the ``latitude--longitude
grid'' \cite{902-Swinbank,948-Williamson} or ``equal-angle grid''
\cite{Gregory2008}. The points concentrate towards the poles, due
to the converging meridians, resulting in high anisotropy.

The number of points, $P$, depends on the angular spacing,
$\delta$, between grid lines. Since $\delta=180^\circ/k$ with
$k=1,2,\dots$,
\begin{equation}
P=2k(k-1)+2.
\end{equation}
That is the number of meridians $(2k)$ times the number of
parallels $(k-1)$, plus the two poles. Frequently, to evaluate
satellite coverage, \cite{938-Feng} $\delta=0.25^{\circ}$, so more
than a million points are used.

\section{Fibonacci Lattice}\label{sec:FiboLattice}

The points of the Fibonacci lattice are arranged along a tightly
wound generative spiral, with each point fitted into the largest
gap between the previous points (Fig.~\ref{fig:GenSpirals}). This
spiral is not apparent (Fig.~\ref{fig:SphereLattices}) because
consecutive points are spaced far apart from each other. The most
apparent spirals join the points closest to each other, and form
crisscrossing sets \cite{902-Swinbank}. The points are evenly
spaced in a very isotropous way.

\begin{figure}[t]
\begin{center}
\includegraphics[width=\textwidth]{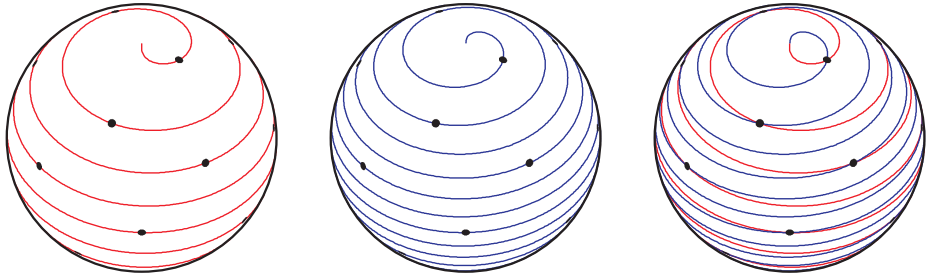}
\end{center}
\textbf{\caption{Generative spirals of a Fibonacci lattice with 21
points\textnormal{. The angle turn between consecutive points
along a spiral is based on the golden ratio ($\phi$):~either the
golden angle, $360^\circ\phi^{-2}\simeq 137.5^\circ$
(\textit{first spiral, red}), or its complementary,
$360^\circ\phi^{-1}\simeq 222.5^\circ$ (\textit{second spiral,
blue}). No point is placed at the poles. Orthographic projection,
centred at longitude $0^\circ$, latitude
$45^\circ$.}}\label{fig:GenSpirals}}
\end{figure}

The next subsections describe how to construct the lattice used in
this paper, and the history of the Fibonacci lattice in various
research fields.

\subsection{Lattice Construction}

The Fibonacci lattice is named after the Fibonacci ratio. The
Fibonacci sequence was discovered in ancient India
\cite{Singh1985,Knuth1997} and rediscovered in the middle ages by
Leonardo Pisano, better known by his nickname Fibonacci
\cite{Sigler2002}. Each term of the sequence, from the third
onwards, is the sum of the previous two: 0, 1, 1, 2, 3, 5, 8, 13,
21, $\dots$ Given two consecutive terms, $F_{i}$ and $F_{i+1}$, a
Fibonacci ratio is $F_{i}/F_{i+1}$. As first proved by Robert
Simson in 1753 \cite{530-Weisstein}, this quotient, as $i\to
\infty$, quickly approaches the golden ratio, defined as
$\phi=1+\phi^{-1}=(1+\sqrt{5})/2\simeq 1.618$.

The Fibonacci lattice differs from other spiral lattices on the
sphere
\cite{Weiller1966,Klima1981,982-Rakhmanov,916-Saff,949-Chukkapalli,981-Bauer,950-Huettig}
in that the longitudinal turn between consecutive points along the
generative spiral is the golden angle,
$360^\circ(1-\phi^{-1})=360^\circ\phi^{-2}\simeq 137.5^{\circ}$,
or its complementary, $360^\circ\phi^{-1}\simeq 222.5^{\circ}$.
Some lattice versions replace $\phi$ by its rational approximant,
a Fibonacci ratio.

The golden angle optimizes the packing efficiency of elements in
spiral lattices \cite{986-Ridley,987-Ridley}. This is because the
golden ratio is the most irrational number \cite{530-Weisstein},
so periodicities or near-periodicities in the spiral arrangement
are avoided, and clumping of the lattice points never occurs
\cite{986-Ridley,Dixon1987,990-Prusinkiewicz,Jean1994,901-Hannay}.

The lattice version used here \cite{902-Swinbank} is probably the
most homogeneous one. It is generated with a Fermat spiral (also
known as the cyclotron spiral)
\cite{969-Vogel,Dixon1987,Dixon1992}, which embraces an equal area
per equal angle turn, so the area between consecutive sampling
points, measured along the spiral, is always the same
\cite{902-Swinbank}. Also, its first and last points are offset
from the poles, leading to a more homogeneous polar arrangement
\cite{970-Purser,902-Swinbank} than in other versions
\cite{974-Svergun,962-Kozin,977-Nye,901-Hannay,970-Purser}. When a
Fibonacci ratio is used
\cite{974-Svergun,962-Kozin,977-Nye,901-Hannay}, the number of
lattice points is $F+1$, where $F>1$ is a term of the Fibonacci
sequence. The lattice used here is instead based on the golden
ratio and can have any odd number of points.

To elaborate the lattice \cite{902-Swinbank}, let $N$ be any
natural number. Let the integer $i$ range from $-N$ to $+N$. The
number of points is
\begin{equation}
P=2N+1,
\end{equation}
and the spherical coordinates, in radians, of the $i$th point are:
\begin{eqnarray}
\mathrm{lat}_i &=& \arcsin \left(\frac{2i}{2N+1} \right),\label{eq:FiboLat}\\
\mathrm{lon}_i &=& 2\pi i\phi^{-1}.\label{eq:FiboComplGolden}
\end{eqnarray}

This pseudocode provides the geographical coordinates in degrees:\\
For $i = -N, (-N+1), \dots, (N-1), N$, Do \{ \\
\indent $\mathrm{lat}_i = \arcsin \left(\frac{2i}{2N+1} \right)\times 180^{\circ}/\pi$\\
\indent $\mathrm{lon}_i = \textrm{mod}(i, \phi) \times 360^\circ/\phi$\\
\indent If $\mathrm{lon}_i < -180^\circ$, then $\mathrm{lon}_i=360^\circ+\mathrm{lon}_i$\\
\indent If $\mathrm{lon}_i > 180 ^\circ$, then $\mathrm{lon}_i=\mathrm{lon}_i-360^\circ$\\
\} End Do\\
Here, arcsin returns a value in radians, while $\textrm{mod}(x,y)$
returns the remainder when $x$ is divided by $y$, removing the
extra turns of the generative spiral. For example, $\textrm{mod}(
6, \phi)=6-3\times \phi$. The last two lines keep the longitude
range from $-180^{\circ}$ to $+180^{\circ}$.

Every point of this lattice is located at a different latitude,
providing a more efficient sampling than the latitude--longitude
lattice. The middle point, $i=0$, is placed at the equator
($\mathrm{lat}_0=0$ and $\mathrm{lon}_0=0$). Each of the other
points $(\mathrm{lat}_i,\mathrm{lon}_i)$ with $i\ne0$, has a
centrosymmetric one with $(-\mathrm{lat}_i,-\mathrm{lon}_i)$. The
lattice as a whole is not centrosymmetric.

The longitudinal turn between consecutive lattice points along the
spiral (Eq.~\ref{eq:FiboComplGolden}) is the complementary of the
golden angle \cite{Swinbank1999,902-Swinbank}. To use the golden
angle instead, we can substitute Eq.~\ref{eq:FiboComplGolden} by
\begin{equation}
\mathrm{lon}_i = -2\pi i\phi^{-2}.\label{eq:FiboGolden}
\end{equation}
With Eq.~\ref{eq:FiboComplGolden}, the spiral progresses
eastwards, while the minus sign of Eq.~\ref{eq:FiboGolden}
indicates a westward progression.

A remark not found elsewhere is that the lattice points are placed
at the intersections between these Fermat spirals of opposite
chirality, except at the poles (Fig.~\ref{fig:GenSpirals}). For
drawing the spirals, $i$ is made continuous in
Eqs.~\ref{eq:FiboLat}, \ref{eq:FiboComplGolden} and
\ref{eq:FiboGolden}, and ranges from $(-N-1/2)$ to $(N+1/2)$. The
$1/2$ term accounts for the polar offset
\cite{970-Purser,902-Swinbank}.

\subsection{Lattice History}\label{sec:FiboReview}

The Fibonacci lattice is a mathematical idealization of patterns
of repeated plant elements, such as rose petals, pineapple scales,
or sunflower seeds (Fig.~\ref{fig:Mammillaria}). The study of
these arrangements is known as Phyllotaxis
\cite{990-Prusinkiewicz,Jean1994,984-Adler,989-Kuhlemeier}. The
Bravais brothers \cite{Bravais1837} were the first to describe
them using a spiral lattice on a cylinder. They argued that the
most common angle turn between consecutive elements along this
spiral in plants is the golden angle. The latter provides optimum
packing \cite{986-Ridley,987-Ridley}, maximizing the exposure to
light, rain and insects for pollination \cite{955-Macia}.
Structures in cells and viruses also follow this pattern
\cite{Erikson1973}. In some experiments, elements are
spontaneously ordered on roughly hemispherical Fibonacci lattices,
because the system tends to minimize the strain energy
\cite{956-Li} or to avoid periodic organization \cite{947-Douday}.

\begin{figure}
\begin{center}
\includegraphics[width=0.7\textwidth]{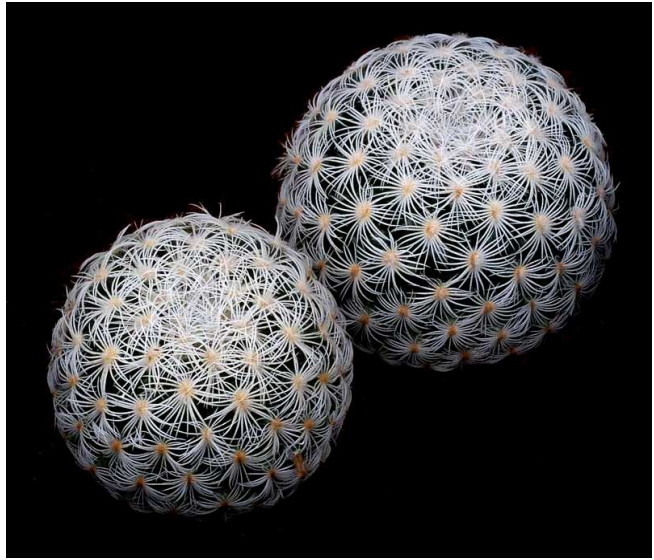}
\end{center}
\textbf{\caption{Spherical Fibonacci lattices in Nature
\textnormal{-- an oblique view of two \emph{Mammillaria
solisioides}. The elements of many plants form Fibonacci lattices,
for example, the areoles (spine-bearing nodes) of
\emph{Mammillaria} cacti
\cite{988-Fowler,989-Kuhlemeier}.}}\label{fig:Mammillaria}}
\end{figure}

\begin{figure}
\begin{center}
\includegraphics[width=0.7\textwidth]{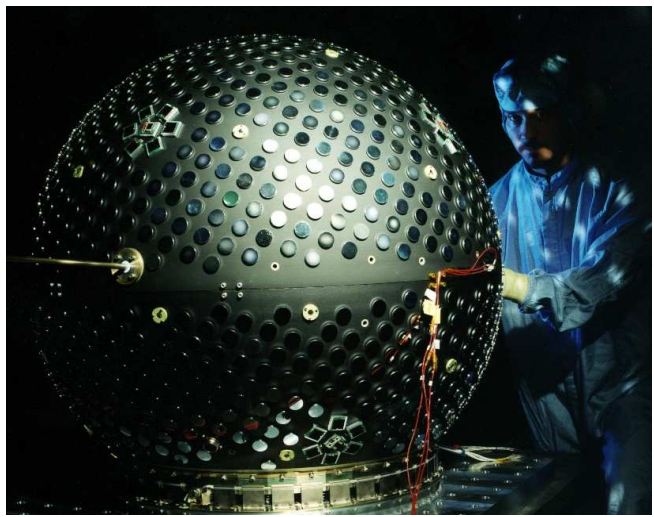}
\end{center}
\textbf{\caption{The Starshine-3 satellite \textnormal{
\cite{980-Maley,991-Lean} had 1500 mirrors arranged on a spherical
Fibonacci lattice. Picture by Michael A. Savell and Gayle R.
Fullerton, taken while the satellite was being inspected by John
Vasquez. Reproduced by courtesy of the U.S. Naval Research
Laboratory.}}\label{fig:Starshine}}
\end{figure}

Unwrapping the cylindrical Fibonacci lattice produces a flat one
\cite{Bravais1837,901-Hannay,902-Swinbank}, frequently used for
numerical integration
\cite{Zaremba1966,Niederreiter1992,Niederreiter1994,Sloan1994}.

Projecting the cylindrical Fibonacci lattice to the sphere
generates the spherical version \cite{901-Hannay,902-Swinbank}.
This can be generalized to arbitrary surfaces of revolution
\cite{987-Ridley,Dixon1992}. The first graphs of spherical
Fibonacci lattices used, as here, the golden ratio and a Fermat
spiral \cite{Dixon1987,Dixon1989,Dixon1992}. A version based on
the Fibonacci ratio is used in the modelling of complex molecules
\cite{985-Vriend,974-Svergun,962-Kozin}. In this case
\cite{974-Svergun,962-Kozin}, a ``+'' sign in the formula for the
longitude should be substituted by ``$\times$'' (D. Svergun,
personal communication, 2009). Versions with the golden ratio
serve to simulate plants realistically \cite{988-Fowler} and to
design golf balls \cite{978-Winfield}. The latter method was used
by Douglas C. Winfield (B. Moore and D. C. Winfield, personal
communication, 2009) in the Starshine-3 satellite
(Fig.~\ref{fig:Starshine}).

The spherical Fibonacci lattice is a highly efficient sampling
scheme for integrating continuous functions
\cite{977-Nye,901-Hannay,970-Purser}, as was observed in magnetic
resonance imaging \cite{946-Ahmad}. It is also advantageous for
providing grid nodes in global meteorological models
\cite{968-Michalakes,Swinbank1999,Swinbank2006a,902-Swinbank,970-Purser,967-Purser}.

\section{Area Measurement}\label{sec:AreaMeasure}

The area measurement starts by finding which points of the lattice
are inside the considered region. This is expressed by a Boolean
function $f_i$, such that $f_i=1$ if the $i$th point is inside,
and $f_i=0$ otherwise \cite{939-Kantsiper}. If the region is a cap
with radius $r$, it suffices to measure $d$, the shortest distance
(great-circle distance) between the sampling point and the cap
centre (Fig.~\ref{fig:CapMeasure}):
\begin{equation}
f_i = \left\{ \begin{array}{ll} 1 & \textrm{if $d\le r$},\\
0 & \textrm{if $d>r$}. \end{array} \right.
\end{equation}

\begin{figure}[!ht]
\begin{center}
\includegraphics[width=0.6\textwidth]{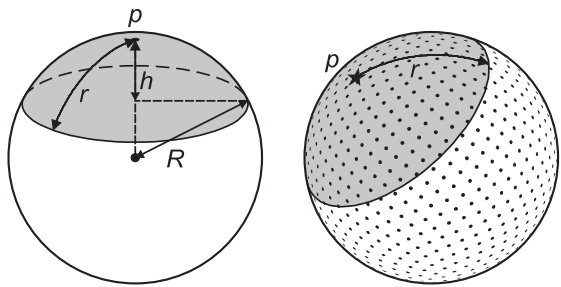}
\end{center}
\textbf{\caption{Measurement of cap area\textnormal{. The cap
(\textit{left}), centred at $p$, has height $h$ and great-circle
radius $r$, and is a region of a sphere with radius $R$. Placed at
random (\textit{right}), its area can be estimated by considering
the lattice points it covers.}}\label{fig:CapMeasure}}
\end{figure}

Each lattice point must be assigned a weight, $w_i$, proportional
to the area it represents. Then the estimate, $\widetilde{A}$, of
the region area $A$, is measured considering the sphere area,
$A_S$, and summing the contribution of all $P$ points of the
lattice:
\begin{equation}\label{eq:AreaEstimate}
A \simeq \widetilde{A}= \frac{A_S \sum_{i=1} ^{P} f_i
w_i}{\sum_{i=1} ^{P} w_i}.
\end{equation}
The weights depend on the lattice type, as described below.

\subsection{Weights in the Latitude--Longitude Lattice}

The weights should be inversely proportional to the point density,
which here increases towards the poles
(Fig.~\ref{fig:SphereLattices}). The linear spacing between
parallels is constant. The length of a parallel is $2\pi R
\cos(\mathrm{lat})$, where $R$ is the sphere radius. In any
parallel there is the same number of lattice points ($2k$), so
their density is inversely proportional to
\cite{939-Kantsiper,973-VandenDool}:
\begin{equation}
w_i=\cos(\mathrm{lat}_i)
\end{equation}

\subsection{Weights in the Fibonacci lattice}

Thanks to the even distribution of points in this lattice, the
same weight can be assumed for all of them
\cite{970-Purser,902-Swinbank}, namely
\begin{equation}
w_i=1.
\end{equation}

Each point represents the area corresponding to its Voronoi cell
(Thiessen polygon). This is the region of positions closer to the
corresponding lattice point than to any other
\cite{Evans1987,Na2002}. Using the exact area of each Voronoi cell
as weight for its lattice point \cite{946-Ahmad} would improve the
area measurement only slightly. The average cell area equals
$A_S/P$. Here, regardless of $P$, only the areas of less than ten
cells, located at the polar regions, differ by more than $\sim
2\%$ from this value \cite{902-Swinbank}.  As $P$ increases,
proportionally fewer cells depart significantly from the average
area. Unlike the latitude--longitude lattice, the homogeneity of
the Fibonacci lattice improves with the number of points.

\section{Error Assessment}\label{sec:ErrorAssess}

This section details how to assess the error involved in measuring
the area of spherical caps placed at random on the sphere. A good
way to measure the homogeneity of a spherical lattice is to
compare the proportion of lattice points in spherical regions with
the normalized areas of the regions \cite{983-Cui}. For this task,
it is natural to use spherical caps
\cite{916-Saff,958-Damelin,Brauchart2004}, which moreover appear
in the applications mentioned in the introduction.

The area of a spherical cap (see Fig.~\ref{fig:CapMeasure}) is:
\begin{equation}
A_C=2\pi R h=2\pi R^2 \left(1- \cos \frac{r}{R} \right).
\end{equation}

The normalized cap area is:
\begin{equation}
F=\frac{A_C}{A_S}=\frac{1-\cos (r/R)}{2},
\end{equation}
where $A_S=4\pi R^2$ is the sphere area.

The absolute error of measuring a single cap is the absolute
difference between the estimated fraction and the actual one:
\begin{equation}\label{eq:AbsError}
E=\left|\widetilde{A}_C/A_S-F \right|.
\end{equation}
This depends on the lattice type, the number of points, and the
size and location of the cap. If $E=0$, the cap gets its fair
share of weighted lattice points. If $A_C=0$, $E=0$. A plane cuts
the sphere into two complementary caps. For any cap, $E$ is the
same as for the complementary cap with area $A'_C=1-A_C$, so it
suffices to consider caps not larger than a hemisphere
($A_C=A_S/2$).

The error is characterized here numerically using a Monte Carlo
method. In each particular realization, a cap is randomly placed.
Every point of the sphere has the same probability of becoming the
cap centre, $p$, with coordinates:
\begin{eqnarray}
\mathrm{lat}_p&=&\frac{180^\circ}{\pi}\arcsin(2X-1)\\
\mathrm{lon}_p&=&360^\circ X-180^\circ.
\end{eqnarray}
Here $X$ is a random number, chosen with uniform probability in
the range $[0,1]$, independently for each equation. The area of
this $j$th cap is estimated, and its corresponding error is $E_j$.

The process is repeated for a total of $n$ independently located
caps of the same size, providing a sample of $n$ values of $E$.
The root mean squared error is \cite{530-Weisstein}:
\begin{equation}
\mathrm{rmse}=\sqrt{\frac{1}{n}\sum_{j=1}^{n}E_j^2}.
\end{equation}

The supremum error in Eq.~\ref{eq:AbsError} for caps of any size,
for a given lattice type and $P$, using $w_i=1$, is the
``spherical cap discrepancy''
\cite{916-Saff,958-Damelin,Brauchart2004}. It cannot be determined
exactly with a Monte Carlo simulation because it might result from
a cap size not used, or a location not sampled. However, it is
unfortunately difficult to compute explicitly \cite{958-Damelin}.
The maximum $E$ measured for the Fibonacci lattice is a lower
bound to its spherical cap discrepancy.

\section{Results}\label{sec:ErrResults}

This section details the maximum errors and root mean squared
errors measured with the Monte Carlo simulation detailed in the
previous section.

Thirteen lattice configurations of each type, from $P \simeq 10^2$
to $P \simeq 10^6$, were analyzed. The chosen values of $P$
increase in logarithmic steps, as accurately as possible ($P
\simeq 10^2, 10^{7/3}, 10^{8/3}, 10^3, \dots$). It is impossible
to use identical $P$ for both lattices because $P$ is odd in the
Fibonacci lattice but even in the latitude--longitude lattice.
Moreover, there are considerably fewer possible values of $P$ in
the latter. For each configuration, 200 different cap sizes were
used: from $A_C=0.0025A_S$ to $A_C=0.5A_S$ in steps of
$0.0025A_S$. To obtain smooth results, $n=60,000$ was chosen for
each cap size and lattice configuration.

In total, 312 million caps were measured (2 lattice types $\times$
13 values of $P$ $\times$ 200 cap sizes $\times$ 60,000 caps).
After optimization, the calculation took 43 days of CPU time using
2.8 GHz, 64-bit processors.

The maximum $E$ measured, for caps of any size and location, is
represented in Fig.~\ref{fig:MaxFiboError}. For the Fibonacci
lattice, it is much lower and decays faster than for the
latitude--longitude lattice, despite the non-uniform weighting of
points of the latter.

\begin{figure}[!ht]
\begin{center}
\includegraphics[width=0.8\textwidth]{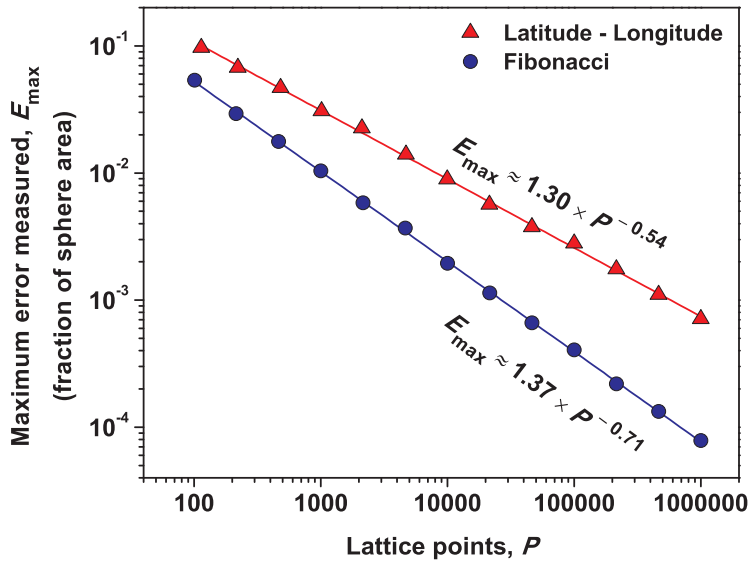}
\end{center}
\textbf{\caption{Maximum error \textnormal{measured for randomly
located spherical caps of any size. For area measurements on the
Earth, about a million points are frequently used \cite{938-Feng},
for which the maximum error would be an order of magnitude smaller
in the Fibonacci lattice.}}\label{fig:MaxFiboError}}
\end{figure}

The rmse depends on the lattice type, number of points and cap
size, as shown in Fig.~\ref{fig:RMSEFibo}. Because of the symmetry
of the latitude--longitude lattice, any hemispherical cap covers
one half of the points of the lattice, and (thanks to the point
weighs) the estimation is perfect ($E=0$ and rmse$=0$). This
exception aside, the rmse tends to increase with the cap area, in
a non-trivial way, which is more complex for the
latitude--longitude lattice than for the Fibonacci lattice.

\begin{figure}
\begin{center}
\includegraphics[width=0.8\textwidth]{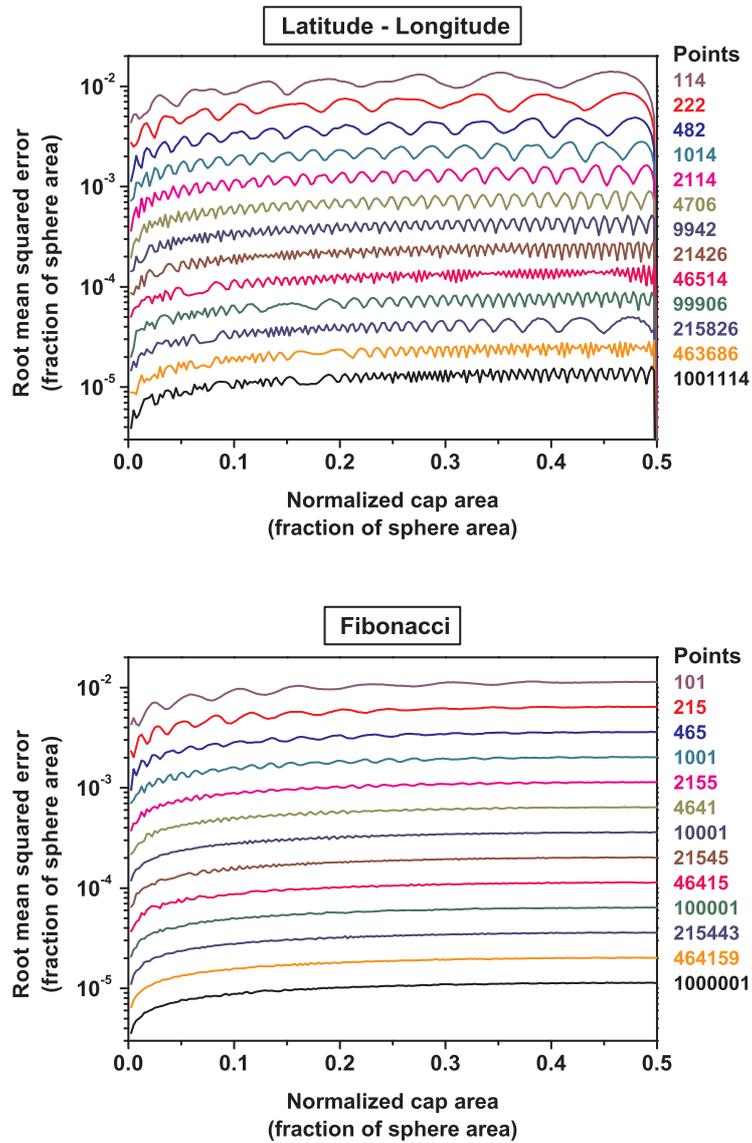}
\end{center}
\textbf{\caption{ Root mean squared error\textnormal{. For
randomly placed caps which occupy an area fraction given by the
abscissas, the curve indicates the root mean squared error of the
area measurement. Each curve corresponds to a different number of
lattice points, labelled at its right. For similar lattice
densities, the Fibonacci lattice provides smaller and more
homogeneous errors.}}\label{fig:RMSEFibo}}
\end{figure}

Figure \ref{fig:MaxRMSEFibo} (top) shows the maximum values of
rmse of each curve. They follow parallel power laws:
\begin{equation}
\mathrm{rmse}_{\mathrm{max}}\simeq k P^{-3/4},
\end{equation}
with $k\simeq 0.505$ for the latitude--longitude lattice, and
$k\simeq 0.362$ for the Fibonacci lattice. Interpolating for the
same $P$, the rmse$_{\mathrm{max}}$ would be $0.505/0.362 \simeq
40\%$ larger in the latitude--longitude lattice.

\begin{figure}[!ht]
\begin{center}
\includegraphics[width=0.8\textwidth]{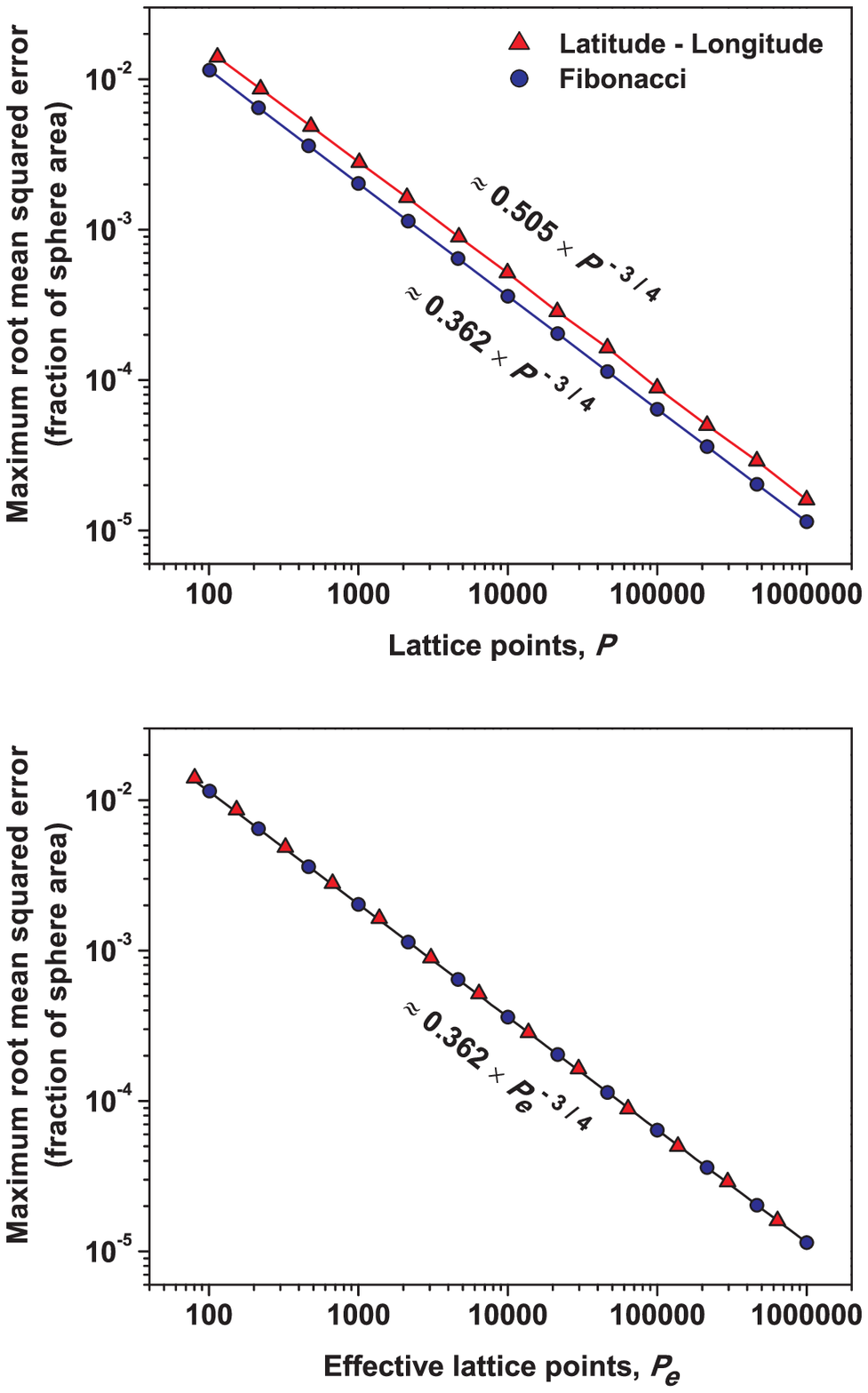}
\end{center}
\textbf{\caption{Scaling of the maximum root mean squared
error\textnormal{.
 Each dot corresponds to the maximum value of a curve
in Fig.~\ref{fig:RMSEFibo}. \textit{Top}: For the same number of
points, the values would be about 40\% smaller for the Fibonacci
lattice. \textit{Bottom}: The latitude--longitude lattice is
inefficient because the effective (weighted) number of points is
smaller than the real one. Considering this fact, the results
collapse to a single function. The power-law decay and its
exponent agree with analytical scaling arguments (see
text).}}\label{fig:MaxRMSEFibo}}
\end{figure}

\subsection{Analytical Approach to the Root Mean Squared Error}

The scaling shown in Fig.~\ref{fig:MaxRMSEFibo} can be explained
using arguments from similar problems in the plane
\cite{Kendall1948,992-Huxley,952-Huxley}. The number of points of
a regular square lattice enclosed by a sufficiently smooth, closed
curve placed at random can be expressed as \cite{952-Huxley}:
\begin{equation}
P_{\mathrm{in}}=AM^2+D,\label{eq:PointsInCurve}
\end{equation}
where $A$ is the area enclosed by the curve, $M$ is the inverse to
the lattice spacing, and $D$ is the discrepancy. The nominal
spacing of a spherical lattice is $\sqrt{A_S/P}$
\cite{902-Swinbank}. Assuming that the Fibonacci lattice is
regular enough,
\begin{equation}
M\approx\sqrt{\frac{P}{A_S}}.
\end{equation}
Substituting this into Eq.~\ref{eq:PointsInCurve},
\begin{equation}
P_{\mathrm{in}}\approx P\frac{A}{A_S}+D.
\end{equation}
Dividing all the terms of this equation by $P$,
\begin{eqnarray}
E \approx \frac{\left|D\right|}{P}\,\,\, ; \,\,\, \mathrm{rmse}
\approx \frac{\mathrm{rms}(D)}{P}.
\end{eqnarray}
In the planar case, the root mean squared of $D$, rms$(D)$, is
proportional to $\sqrt{M}$
\cite{Kendall1948,992-Huxley,952-Huxley}. Extrapolating this fact,
we obtain the scaling observed in the Monte Carlo simulation:
\begin{equation}
\mathrm{rmse} \propto \frac{\sqrt{M}}{P} \approx
\frac{(P/A_S)^{1/4}}{P}\propto P^{-3/4}.
\end{equation}

In the latitude--longitude lattice, the same scaling holds if we
consider its smaller sampling efficiency. The latter may be
measured with the denominator of Eq.~\ref{eq:AreaEstimate}. This
is the number of points of a homogeneous lattice that would do the
same work, and is up to $\sim36\%$ smaller than $P$ in the range
of $P$ considered here. If Fig.~\ref{fig:MaxRMSEFibo} is plotted
using these effective points in the abscissas, the same fit
suffices for both lattice types (Fig.~\ref{fig:MaxRMSEFibo},
bottom).

If the sampling points were placed at random, the rmse would
decrease more slowly, as $\propto P^{-1/2}$ \cite{Bevington1992}.

\section{Conclusions}

This paper analyzes the errors involved in measuring the areas of
spherical caps using lattices of sampling points: the
latitude--longitude lattice (classically used for this task), and
a Fibonacci lattice \cite{902-Swinbank}. The latter has low
anisotropy, is easy to construct, and is shown to result from the
intersection of two generative spirals
(Fig.~\ref{fig:GenSpirals}). A a review of the literature
(Sect.~\ref{sec:FiboReview}) reveals successful applications of
this spherical lattice since the 1980s.

If the Fibonacci lattice were used instead of its
latitude--longitude counterpart, the area measurement would be
more efficient, allowing a significant reduction of the
computation time. For approximately the same number of lattice
points, the maximum root mean squared error would be about 40\%
smaller (Fig.~\ref{fig:MaxRMSEFibo}). The maximum errors would be
also smaller, and would decay faster with the number of points
(Fig.~\ref{fig:MaxFiboError}). If about a million points were
used, as is commonly the case \cite{938-Feng}, the maximum error
would be an order of magnitude smaller
(Fig.~\ref{fig:MaxFiboError}).

It is also found that the maximum root mean squared error obeys a
single scaling relation when the sampling efficiency is taken into
account (Fig.~\ref{fig:MaxRMSEFibo}). This is partially explained
using arguments from similar problems on the plane
\cite{Kendall1948,992-Huxley,952-Huxley}.

The area estimate depends also on the orientation of the sampling
lattice, especially if the latter has high anisotropy. The
difference may be assessed \cite{901-Hannay} by rotating the
lattice \cite{1007-Greiner}. Such an issue is not relevant in the
case analyzed here because the caps were placed at random with
uniform probability over the spherical surface.

Here the Earth's shape has been approximated by a sphere
\cite{939-Kantsiper}, adding a slightly higher error than other
shape models \cite{971-Kimerling,953-Earle,917-Sjoberg}. Assessing
this difference may be a topic of future research.

\section{Acknowledgements}
I thank John H. Hannay, Gil Moore,
Dimitri Svergun, Richard Swinbank, John Vasquez, Gert Vriend and
Douglas C. Winfield for clarifying different aspects of their
work. The results were calculated using the supercomputing
facilities of the Institute of Biocomputing and Physics of Complex
Systems (BIFI, University of Zaragoza, Spain). Lattice maps were
produced with Generic Mapping Tools \cite{Wessel1998}. I am also
grateful to an anonymous reviewer for offering insightful and
encouraging comments.


\begin{thebibliography}{}

\bibitem{984-Adler}
Adler, I.; Barabe, D. \& Jean, R.~V. (1997) A history of the study
of phyllotaxis. \textit{Annals of Botany}, 80 (3), 231--244.
Doi:10.1006/anbo.1997.0422

\bibitem{946-Ahmad}
Ahmad, R.; Deng, Y.; Vikram, D.~S.; Clymer, B.; Srinivasan, P.;
Zweier, J.~L. \& Kuppusamy, P. (2007) Quasi Monte Carlo-based
isotropic distribution of gradient directions for improved
reconstruction quality of 3D EPR imaging. \textit{Journal of
Magnetic Resonance}, 184 (2), 236--245.
Doi:10.1016/j.jmr.2006.10.008

\bibitem{Baddeley2004}
Baddeley, A. \& Jensen, E.~B.~V. (2004) \textit{Stereology for
Statisticians}. CRC Press, Boca Raton, Florida. 416 pp.

\bibitem{972-Barclay}
Barclay, M. \& Galton, A. (2008) Comparison of region
approximation techniques based on Delaunay triangulations and
Voronoi diagrams. \textit{Computers, Environment and Urban
Systems}, 32 (4), 261--267.
Doi:10.1016/j.compenvurbsys.2008.06.003

\bibitem{995-Bardsley}
Bardsley, W.~E. (1983) Random error in point counting.
\textit{Mathematical Geology}, 15 (3), 469--475.
Doi:10.1007/BF01031293

\bibitem{981-Bauer}
Bauer, R. (2000) Distribution of points on a sphere with
application to star catalogs. \textit{Journal of Guidance, Control
and Dynamics}, 23 (1), 130--137.

\bibitem{Bevington1992}
Bevington, P.~R. \& Robinson, D.~K. (1992) \textit{Data Reduction
and Error Analysis for the Physical Sciences}, 2nd edn. McGraw
Hill, Boston, 328 pp.

\bibitem{940-Bevis}
Bevis, M. \& Cambareri, G. (1987) Computing the area of a
spherical polygon of arbitrary shape. \textit{Mathematical
Geology}, 19 (4), 335--346. Doi:10.1007/BF00897843

\bibitem{Brauchart2004}
Brauchart, J.~S. (2004) Invariance principles for energy
functionals on spheres. \textit{Monatshefte f\"ur Mathematik}, 141
(2), 101--117. Doi:10.1007/s00605-002-0007-0

\bibitem{Bravais1837}
Bravais, L. \& Bravais, A. (1837) Essai sur la disposition des
feuilles curvis\'eri\'ees. \textit{Annales des Sciences
Naturelles} 7, 42--110 \& plates 2--3.

\bibitem{949-Chukkapalli}
Chukkapalli, G.; Karpik, S.~R. \& Ethier, C.~R. (1999) A scheme
for generating unstructured grids on spheres with application to
parallel computation. \textit{Journal of Computational Physics},
149 (1), 114--127. Doi:10.1006/jcph.1998.6146

\bibitem{Conway1998}
Conway, J.~H. \& Sloane, N.~J.~A. (1998) \textit{Sphere Packings,
Lattices and Groups}, 3rd edn. Springer-Verlag, New York, 703 pp.

\bibitem{983-Cui}
Cui, J. \& Freeden, W. (1997) Equidistribution on the sphere.
\textit{SIAM Journal on Scientific Computing}, 18 (2), 595--609.
Doi:10.1137/S1064827595281344

\bibitem{958-Damelin}
Damelin, S.~B. \& Grabner, P.~J. (2003) Energy functionals,
numerical integration and asymptotic equidistribution on the
sphere. \textit{Journal of Complexity}, 19 (3), 231--246.
Doi:10.1016/S0885-064X(02)00006-7 [Corrigendum in 20 (6),
883--884. Doi:10.1016/j.jco.2004.06.003]

\bibitem{Dixon1987}
Dixon, R. (1987) \textit{Mathographics}. Basil Blackwell, Oxford,
England, 224 pp.

\bibitem{Dixon1989}
Dixon, R. (1989) Spiral phyllotaxis. \textit{Computers and
Mathematics with Applications}, 17 (4--6), 535--538.

\bibitem{Dixon1992}
Dixon, R. (1992) Green spirals. In:~Hargittai, I. \& Pickover,
C.~A. (Eds.) \textit{Spiral Symmetry}. World Scientific,
Singapore, pp.~353--368.

\bibitem{947-Douday}
Douady, S. \& Couder, Y. (1992) Phyllotaxis as a physical
self-organized growth process. \textit{Physical Review Letters} 68
(13), 2098--2101. Doi:10.1103/PhysRevLett.68.2098

\bibitem{953-Earle}
Earle, M.~A. (2006) Sphere to spheroid comparisons. \textit{The
Journal of Navigation} 59 (3), 491--496.
Doi:10.1017/S0373463306003845

\bibitem{Erikson1973}
Erikson, R.~O. (1973) Tubular packing of spheres in biological
fine structure. \textit{Science}, 181 (4101), 705--716.
Doi:10.1126/science.181.4101.705

\bibitem{Evans1987}
Evans, D.~G. \& Jones, S.~M. (1987) Detecting Voronoi
(area-of-influence) polygons. \textit{Mathematical Geology}, 19
(6), 523--537. Doi:10.1007/BF00896918

\bibitem{938-Feng}
Feng, S.; Ochieng, W.~Y. \& Mautz, R. (2006) An area computation
based method for RAIM holes assessment. \textit{Journal of Global
Positioning Systems}, 5 (1--2), 11--16.

\bibitem{988-Fowler}
Fowler, D.~R.; Prusinkiewicz, P. \& Battjes, J. (1992) A
collision-based model of spiral phyllotaxis. \textit{ACM SIGGRAPH
Computer Graphics}, 26 (2), 361--368. Doi:10.1145/142920.134093

\bibitem{GonzalezStatsei6}
Gonz\'alez, \'A. (2009) Self-sharpening seismicity maps for
forecasting earthquake locations. In:~\textit{Abstracts of the
Sixth International Workshop on Statistical Seismology, Tahoe
City, California, 16--19 April 2009}.
\href{http://www.scec.org/statsei6/posters.html}
{http://www.scec.org/statsei6/posters.html} [Last accessed:
November 16th, 2009]

\bibitem{Gregory2008}
Gregory, M.~J.; Kimerling, A.~J; White, D. \& Sahr, K. (2008) A
comparison of intercell metrics on discrete global grid systems.
\textit{Computers, Environment and Urban Systems}, 32 (3),
188--203. Doi:10.1016/j.compenvurbsys.2007.11.003

\bibitem{1007-Greiner}
Greiner, B. (1999) Euler rotations in plate-tectonic
reconstructions. \textit{Computers and Geosciences}, 25 (3),
209--216. Doi:10.1016/S0098-3004(98)00160-5

\bibitem{997-Gundersen}
Gundersen, H.~J.~G.; Jensen, E.~B.~V.; Ki\^eu, K. \& Nielsen, J.
(1999) The efficiency of systematic sampling in stereology --
reconsidered. \textit{Journal of Microscopy}, 193 (3), 199-211.
Doi:10.1046/j.1365-2818.1999.00457.x

\bibitem{901-Hannay}
Hannay, J.~H. \& Nye, J.~F. (2004) Fibonacci numerical integration
on a sphere. \textit{Journal of Physics A: Mathematical and
General}, 37 (48), 11591--11601. Doi:10.1088/0305-4470/37/48/005

\bibitem{998-Howarth}
Howarth, R.~J. (1998) Improved estimators of uncertainty in
proportions, point-counting, and pass-fail test results.
\textit{American Journal of Science}, 298 (7), 594--607.

\bibitem{950-Huettig}
H\"uttig, C. \& Stemmer, K. (2008) The spiral grid:~a new approach
to discretize the sphere and its application to mantle convection.
\textit{Geochemistry, Geophysics, Geosystems} 9 (2), Q02018.
Doi:10.1029/2007GC001581

\bibitem{992-Huxley}
Huxley, M.~N. (1987) The area within a curve. \textit{Proceedings
of the Indian Academy of Sciences (Mathematical Sciences)}, 97
(1--3), 111--116. Doi:10.1007/BF02837818

\bibitem{952-Huxley}
Huxley, M.~N. (2003) Exponential sums and lattice points III.
\textit{Proceedings of the London Mathematical Society}, 87 (3),
591--609. Doi:10.1112/S0024611503014485

\bibitem{996-Jarai}
Jar\'ai, A.; Koz\'ak, M. \& R\'ozsa, P. (1997) Comparison of the
methods of rock-microscopic grain-size determination and
quantitative analysis. \textit{Mathematical Geology}, 29 (8),
977--991. Doi:10.1023/A:1022305518696

\bibitem{Jean1994}
Jean, R.~V. (1994) \textit{Phyllotaxis:~A Systemic Study of Plant
Pattern Morphogenesis}. Cambridge University Press, Cambridge, UK,
400 pp.

\bibitem{943-Kafka}
Kafka, A.~L. (2007) Does seismicity delineate zones where future
large earthquakes are likely to occur in intraplate environments?
In:~Stein, S. \& Mazzotti, S. (Eds.) \textit{Continental
Intraplate Earthquakes: Science, Hazard, and Policy Issues}.
Geological Society of America Special Paper 425, Boulder,
Colorado, pp.~35--48. Doi:10.1130/2007.2425(03)

\bibitem{939-Kantsiper}
Kantsiper, B. \& Weiss, S. (1998) An analytic approach to
calculating Earth coverage. \textit{Advances in the Astronautical
Sciences}, 97, 313--332.

\bibitem{Kendall1948}
Kendall, D.~G. (1948) On the number of lattice points inside a
random oval. \textit{The Quarterly Journal of Mathematics
(Oxford)}, 19 (1), 1--26. Doi:10.1093/qmath/os-19.1.1

\bibitem{971-Kimerling}
Kimerling, A.~J. (1984) Area computation from geodetic coordinates
on the spheroid. \textit{Surveying and Mapping}, 44 (4), 343--351.

\bibitem{Klima1981}
Kl\'ima, K.; Pick, M. \& Pros, Z. (1981) On the problem of equal
area block on a sphere. \textit{Studia Geophysica et Geodaetica
(Praha)}, 25 (1), 24--35. Doi:10.1007/BF01613559

\bibitem{Knuth1997}
Knuth, D.~E. (1997) \textit{Art of Computer Programming, Volume 1:
Fundamental Algorithms}, 3rd edn. Addison-Wesley, Reading,
Massachusetts, 672 pp.

\bibitem{Kossobokov2003}
Kossobokov, V. \& Shebalin, P. (2003) Earthquake prediction.
In:~Keilis-Borok, V.~I. \& Soloviev, A.~A. (Eds.)
\textit{Nonlinear Dynamics of the Lithosphere and Earthquake
Prediction}. Springer, Berlin, pp.~141--207. [References in
pp.~311--332]

\bibitem{962-Kozin}
Kozin, M.~B.; Volkov, V.~V. \& Svergun, D.~I. (1997) ASSA, a
program for three-dimensional rendering in solution scattering
from biopolymers. \textit{Journal of Applied Crystallography}, 30
(5), 811--815. Doi:10.1107/S0021889897001830

\bibitem{989-Kuhlemeier}
Kuhlemeier, C. (2007) Phyllotaxis. \textit{Trends in Plant
Science}, 12 (4), 143--150. Doi:10.1016/j.tplants.2007.03.004

\bibitem{991-Lean}
Lean, J.~L.; Picone, J.~M.; Emmert, J.~T. \& Moore, G. (2006)
Thermospheric densities derived from spacecraft
orbits:~application to the Starshine satellites. \textit{Journal
of Geophysical Research}, 111 (4), A04301.
Doi:10.1029/2005JA011399

\bibitem{956-Li}
Li, C.; Zhang, X. \& Cao, Z. (2005) Triangular and Fibonacci
number patterns driven by stress on core/shell microstructures.
\textit{Science}, 309 (5736), 909--911.
Doi:10.1126/science.1113412

\bibitem{955-Macia}
Maci\'a, E. (2006) The role of aperiodic order in science and
technology. \textit{Reports on Progress in Physics}, 69 (2),
397--441. Doi:10.1088/0034-4885/69/2/R03

\bibitem{980-Maley}
Maley, P.~D.; Moore, R.~G. \& King, D.~J. (2002) Starshine:~a
student-tracked atmospheric research satellite. \textit{Acta
Astronautica}, 51 (10), 715--721.
Doi:10.1016/S0094-5765(02)00021-8

\bibitem{968-Michalakes}
Michalakes, J.~G.; Purser, R.~J. \& Swinbank, R. (1999) Data
structure and parallel decomposition considerations on a Fibonacci
grid. In:~\textit{Preprints of the 13th Conference on Numerical
Weather Prediction, Denver, 13-17 September 1999}. American
Meteorological Society, pp.~129--130.

\bibitem{Na2002}
Na, H.~S.; Lee, C.~N. \& Cheong, O. (2002) Voronoi diagrams on the
sphere. \textit{Computational Geometry -- Theory and
Applications}, 23 (2), 183--194. Doi:10.1016/S0925-7721(02)00077-9

\bibitem{Niederreiter1992}
Niederreiter, H. (1992) \textit{Random Number Generation and
Quasi-Monte Carlo Methods}. Society for Industrial and Applied
Mathematics, Philadelphia, 247 pp.

\bibitem{Niederreiter1994}
Niederreiter, H. \& Sloan, I.~H. (1994) Integration of nonperiodic
functions of two variables by Fibonacci lattice rules.
\textit{Journal of Computational and Applied Mathematics}, 51 (1),
57--70. Doi:10.1016/0377-0427(92)00004-S

\bibitem{977-Nye}
Nye, J.~F. (2003) A simple method of spherical near-field scanning
to measure the far fields of antennas or passive scatterers.
\textit{IEEE Transactions on Antennas and Propagation}, 51 (8),
2091--2098. Doi:10.1109/TSP.2003.815442

\bibitem{937-Ochieng}
Ochieng, W.~Y.; Sheridan, K.~F.; Sauer, K. \& Han, X. (2002) An
assessment of the RAIM performance of a combined Galileo/GPS
navigation system using the marginally detectable errors (MDE)
algorithm. \textit{GPS Solutions} 5 (3), 42--51.
Doi:10.1007/PL00012898

\bibitem{990-Prusinkiewicz}
Prusinkiewicz, P. \& Lindenmayer, A. (1990) \textit{The
Algorithmic Beauty of Plants}. Springer-Verlag, New York, 228 pp.

\bibitem{967-Purser}
Purser, R.~J. (2008) \textit{Generalized Fibonacci grids:~a new
class of structured, smoothly adaptive multi-dimensional
computational lattices}. Office Note 455, National Centers for
Environmental Prediction, Camp Springs, Maryland, USA, 37 pp.

\bibitem{970-Purser}
Purser, R.~J. \& Swinbank, R. (2006) \textit{Generalized
Euler-Maclaurin formulae and end corrections for accurate
quadrature on Fibonacci grids}. Office Note 448, National Centers
for Environmental Prediction, Camp Springs, Maryland, USA, 19 pp.

\bibitem{982-Rakhmanov}
Rakhmanov, E.~A.; Saff, E.~B. \& Zhou, Y.~M. (1994) Minimal
discrete energy on the sphere. \textit{Mathematical Research
Letters}, 1 (6), 647--662.

\bibitem{986-Ridley}
Ridley, J.~N. (1982) Packing efficiency in sunflower heads.
\textit{Mathematical Biosciences}, 58 (1), 129--139.
Doi:10.1016/0025-5564(82)90056-6

\bibitem{987-Ridley}
Ridley, J.~N. (1986) Ideal phyllotaxis on general surfaces of
revolution. \textit{Mathematical Biosciences}, 79 (1), 1--24.
Doi:10.1016/0025-5564(86)90013-1

\bibitem{916-Saff}
Saff, E.~B. \& Kuijlaars, A.~B.~J. (1997) Distributing many points
on a sphere. \textit{The Mathematical Intelligencer}, 19 (1),
5--11. Doi:10.1007/BF03024331

\bibitem{Sigler2002}
Sigler, L.~E. (2002) \textit{Fibonacci's Liber Abaci:~a
Translation into Modern English of Leonardo Pisano's Book of
Calculation}. Springer, New York, 636 pp.

\bibitem{Singh1985}
Singh, P. (1985) The so-called Fibonacci numbers in ancient and
medieval India. \textit{Historia Mathematica}, 12 (3), 229--244.
Doi:10.1016/0315-0860(85)90021-7

\bibitem{917-Sjoberg}
Sj\"oberg, L.~E. (2006) Determination of areas on the plane,
sphere and ellipsoid. \textit{Survey Review}, 38 (301), 583--593.

\bibitem{Sloan1994}
Sloan, I.~H. \& Joe, S. (1994) \textit{Lattice Methods for
Multiple Integration}. Oxford University Press, Oxford, UK, 256
pp.

\bibitem{974-Svergun}
Svergun, D.~I. (1994) Solution scattering from biopolymers:
advanced contrast-variation data analysis. \textit{Acta
Crystallographica}, A50 (3), 391--402.
Doi:10.1107/S0108767393013492

\bibitem{Swinbank1999}
Swinbank, R. \& Purser, R.~J. (1999) Fibonacci grids.
In:~\textit{Preprints of the 13th Conference on Numerical Weather
Prediction, Denver 13--17 September 1999}. American Meteorological
Society, pp. 125--128.

\bibitem{Swinbank2006a}
Swinbank, R. \& Purser, R.~J. (2006a) \textit{Standard test
results for a shallow water equation model on the Fibonacci grid}.
Forecasting Research Technical Report 480, Met Office, Exeter, UK,
12 pp.

\bibitem{902-Swinbank}
Swinbank, R. \& Purser, R.~J. (2006b) Fibonacci grids:~a novel
approach to global modelling. \textit{Quarterly Journal of the
Royal Meteorological Society}, 132 (619), 1769--1793.
Doi:10.1256/qj.05.227

\bibitem{973-VandenDool}
Van den Dool, H.~M. (2007) \textit{Empirical Methods in Short-Term
Climate Prediction}. Oxford University Press, Oxford, UK, 240 pp.

\bibitem{969-Vogel}
Vogel, H. (1979) A better way to construct the sunflower head.
\textit{Mathematical Biosciences}, 44 (3--4), 179--189.
Doi:10.1016/0025-5564(79)90080-4

\bibitem{985-Vriend}
Vriend, G. (1990) WHAT IF:~a molecular modeling and drug design
program. \textit{Journal of Molecular Graphics}, 8 (1), 52--56.
Doi:10.1016/0263-7855(90)80070-V

\bibitem{Weiller1966}
Weiller, A.~R. (1966) Probleme de l'implantation d'une grille sur
une sphere, deuxi\`eme partie. \textit{Bulletin G\'eod\'esique},
80 (1), 99--111. Doi:10.1007/BF02527041

\bibitem{530-Weisstein}
Weisstein, E.~W. (2002) \textit{CRC Concise Encyclopedia of
Mathematics CD-ROM}, 2nd edn. CRC Press, Boca Raton, Florida, 3252
pp.

\bibitem{Wessel1998}
Wessel, P. \& Smith, W.~H.~F. (1998) New, improved version of
Generic Mapping Tools released. \textit{Eos, Transactions,
American Geophysical Union}, 79 (47), 579.

\bibitem{948-Williamson}
Williamson, D.~L. (2007) The evolution of dynamical cores for
global atmospheric models. \textit{Journal of the Meteorological
Society of Japan}, 85B, 241--269.

\bibitem{978-Winfield}
Winfield, D.~C. \& Harris, K.~M. (2001) \textit{Phyllotaxis-Based
Dimple Patterns}. Patent number WO 01/26749 Al. World Intellectual
Property Organization, Geneva, Switzerland.

\bibitem{Zaremba1966}
Zaremba, S.~K. (1966) Good lattice points, discrepancy, and
numerical integration. \textit{Annali di Matematica Pura ed
Applicata}, 73 (1), 293--317. Doi:10.1007/BF02415091

\end{thebibliography}
\end{document}